\documentclass{article}
\usepackage[utf8]{inputenc}
\usepackage{amsfonts, amsmath, amssymb, amsthm, authblk, graphicx, enumitem}
\usepackage{cite}
\usepackage[colorlinks=true,urlcolor=blue, citecolor=red,linkcolor=blue]{hyperref}
\usepackage[left=2cm, right=2cm, top=3cm]{geometry}
\usepackage[capitalise]{cleveref}
\newtheorem{theorem}{Theorem}

\newtheorem{definition}[theorem]{Definition}

\newtheorem{lemma}[theorem]{Lemma}
\newtheorem{proposition}[theorem]{Proposition}
\newtheorem{remark}[theorem]{Remark}
\usepackage{dsfont}
\usepackage{hyperref}

\renewcommand{\d}{\displaystyle}
\usepackage{graphicx}
\usepackage{apptools}
\AtAppendix{\counterwithin{lem}{section}}

\newcommand{\pts}[1]{\left(#1\right)}  
\newcommand{\cts}[1]{\left[#1\right]}                                  	  
\newcommand{\lvs}[1]{\left\{#1\right\}}                                	  
\newcommand{\abs}[1]{\left|#1\right|}                                  	  

\newcommand{\C}{\mathbb{C}} 
\newcommand{\R}{\mathbb{R}}
\newcommand{\A}{\mathcal{A}}
\newcommand{\al}{\alpha}

\newcommand{\OO}{\mathcal{O}}
\newcommand{\Fi}{\Phi}

\newcommand{\HH}{\mathcal{H}}

\begin{document}
\title{Boundary controllability for a fourth order degenerate parabolic equation with a singular potential 
}

\author{Leandro Galo-Mendoza\thanks{jesus.galo@im.unam.mx}}
\affil{Unidad Cuernavaca, Instituto de Matem\'aticas, Universidad Nacional Aut\'onoma de
M\'exico, M\'exico.}

\maketitle

\abstract{In this paper, we prove the null controllability of a one-dimensional fourth-order degenerate parabolic equation with a singular potential. We analyze the cases when boundary control conditions are applied at the left endpoint. We use a spectral decomposition involving Bessel functions and their zeros in a suitably  weighted Sobolev space for a degenerate parabolic operator with certain boundary conditions. We prove the well-posedness of the system using semigroup operator theory. Then, we apply the moment method by Fattorini and Russell to obtain an upper estimate of the cost of controllability. Additionally, we also obtain a lower estimate of the cost of controllability using a representation theorem for analytic functions of exponential type.}

\section{Introduction and main results}
The aim of this paper is to analyze the null controllability of a certain one-dimensional fourth-order degenerate parabolic equation with a singular potential. Several works have studied the controllability of a one-dimensional fourth-order partial differential equations of parabolic type without degeneracy, for example, see \cite{Carreño}, \cite{MLopez}, \cite{SGuerrero} and \cite{GuoYung}. Some of these works investigate null controllability using more than one boundary control condition or a distributed control condition. Additionally, in these studies, the cost of controllability is examined using Carleman estimates (see \cite{Carreño} and \cite{SGuerrero}) or other techniques (see \cite{MLopez} and \cite{Yu}). Upon reviewing the literature, we found that only a few works address the controllability of such systems using a single control, as seen in \cite{MLopez}. When considering a fourth-order one-dimensional degenerate parabolic equation, the existence and uniqueness of solutions have been analyzed in many studies (see \cite{LGiacomelli}, \cite{Tepoyan}, and more recently \cite{Tepoyan2}). However, we did not find any works specifically addressing the controllability of fourth-order degenerate partial differential equations. In this regard, this paper studies the controllability of a fourth-order degenerate parabolic equation with a single control acting at the degenerate/singular left endpoint. This work is motivated by the formula (\ref{partes}) and the previous studies \cite{GaloLopez}, \cite{GaloLopez2}, and \cite{GaloLopez3}.\\

Let $T > 0$ and set $Q := (0, 1) \times (0, T )$. For $\al,\beta\in \R$ with $0\leq\al<2$, 
consider the operator given by 
\begin{equation}\label{opeA}
	\A u:=-(x^\al u_x)_x-\beta x^{\al -1}u_x-\frac{\mu}{x^{2-\al}} u,
\end{equation}
provided that $\mu\in\R$ satisfies
\begin{equation}\label{mucon}
	-\infty<\mu<\mu(\al+\beta),\quad \text{ and }\quad\mu(\delta):=\frac{(1-\delta)^2}{4}, \quad \delta\in\R, 
\end{equation}
where $u$ is given in a convenient dense domain $D(\A)$. This operator has been studied in ~\cite{GaloLopez, GaloLopez2, GaloLopez3} for solving controllability problems with a control acting by means of a left boundary condition at the degenerate/singular endpoint $x = 0$ or at the regular endpoint $x =1$. 

From (\ref{opeA}) we deduce that if $\A u\in D(\A)$, then
\begin{equation}\label{A_square}
	\A^2u:=x^{2\al}u_{xxxx}+\rho_{1}\,x^{2\al-1}u_{xxx}+\rho_{2}\,x^{2\al-2}u_{xx}+\rho_{3}\,x^{2\al-3}u_{x}+\dfrac{\rho_{4}}{x^{4-2\al}}u,
\end{equation}
is well defined, where
\begin{equation}\label{A2consts}
	\begin{array}{lcl}
		\rho_{1}:=\rho_{1}(\al,\beta) &=& 4\al+2\beta,\\
		\rho_{2}:= \rho_{2}(\al,\beta,\mu) &=& (2\al+\beta)(2\al+\beta-1)+(\al+\beta)(\al-1)+2\mu,\\
		\rho_{3}:=\rho_{3}(\al,\beta,\mu) &=& \cts{(\al+\beta)(\al-1)+2\mu}(2\al+\beta-2),\\
		\rho_{4}:= \rho_{4}(\al,\beta,\mu) &=& \mu\cts{(\al-2)(2\al+\beta-3)+\mu}.
	\end{array}
\end{equation}
The diagonalization of this operator appears as a fourth order differential equation of the Bessel type in \cite{kamke}. In \cite{everitt} a special case of (\ref{A_square}) is studied using a different technique from this work, where the authors demonstrate the diagonalization of their operator over the interval $(0,\infty)$.

Naturally, we are interested in analyzing the null controllability of the system 
\begin{equation}
	u_{t} + \A^2 u = 0 \quad \text{ in }\quad Q
\end{equation}
with suitable boundary conditions. We also need to distinguish between the cases $\al + \beta < 1$ and $\al + \beta > 1$ (The case $\al+\beta = 1$ will be treated separately in Section. \ref{sec5}). For this purpose, if $\al,\beta \in\R$ such that $\al+\beta\neq 1$, we define $\ell =\ell(\al,\beta)$ as follows
\begin{equation}\label{index_l}
	\ell:=\ell(\al,\beta) =\left\{\begin{array}{cc} 0 &\text{if } \al+\beta < 1,\\ 1 & \text{if } \al+\beta >1.\end{array}\right. 
\end{equation}  
Firstly, we introduce suitable weighted boundary operators named $C_{\ell}$ and $D_{\ell}$, acting in a certain smooth function $u$ as follows,
\begin{equation}\label{Dirischlet}
	(C_{\ell}\,u)(t):=\lim_{x\to 0^+}x^{\ell-\gamma}\partial^{\ell}_{x} u(x,t)\,\text{ and } (D_{\ell}\,u)(t): =\lim_{x\to 0^+}x^{\ell(\al+\beta-1)/2}u(x,t),  
\end{equation}
where 
\begin{equation}\label{gamma}
	\gamma=\gamma(\alpha, \beta,\mu):=\frac{1-\al-\beta}{2}-\sqrt{\mu(\alpha+\beta)-\mu},
\end{equation}
and $\ell$ is given by (\ref{index_l}).\\

So, for $r =0, 1$ as a general form, we consider the system 
\begin{equation}\label{problem}
	\left\{\begin{aligned}
		u_t+x^{2\al}u_{xxxx}+\rho_{1}\,x^{2\al-1}u_{xxx}+\rho_{2}\,x^{2\al-2}u_{xx}+\rho_{3}\,x^{2\al-3}u_{x}+\dfrac{\rho_{4}}{x^{4-2\al}}u&=0, & & \text { in }Q, \\
		\cts{C_{\ell}\pts{(1-r)u + r\A u}}(t)= f(t),\quad\quad\quad\quad\quad\,\, u(1, t)&=0, & & \text { on }(0, T), \\
		\cts{D_{\ell}\pts{r\,u+(1-r)\A u}}(t) = 0,\quad\quad \pts{x^{\al+\beta}u_x}_x(1,t)&=0, & & \text { on }(0, T), \\
		u(x, 0) &=u_{0}(x), & & \text { on }(0, 1),
	\end{aligned}\right.
\end{equation}
where $\ell$ is given by (\ref{index_l}) and $\rho_1,\rho_2,\rho_3,\rho_4$ are constants given by (\ref{A2consts}).\\

Referring to the last system, for example, when $\al + \beta > 1$ and $r = 0$, then we have $\ell=1$ and this gives us a pair of weighted boundary conditions at the left endpoint associated with the system. The first condition would come to be $\pts{x^{1-\gamma}u_{x}}(0, t) = f(t)$ where we have the control, and the second condition would be $\pts{x^{(\al+\beta-1)/2}\A u}(0, t) = 0$, a necessary boundary condition to solve the problem (\ref{problem}).\\
Thus, we have four possible systems to be solved, depending of the conditions $\al +\beta < 1$ or $\al +\beta >1 $, $r \in \{0,1\}$,  and in consequence of our weighted boundary operators $C_{\ell}$ and $D_{\ell}$.\\

As a first goal for the system (\ref{problem}) in all the cases to be considered, we define a notion of weak solution and demonstrate the well-posedness of the problem in appropriate interpolation spaces for all possible boundary conditions. Note that, in all cases, we have a weighted boundary condition at the degenerate/singular endpoint where the control is acting. We then apply the moment method introduced by Fattorini and Russell in \cite{Fattorini} to prove the null controllability of the system (\ref{problem}) in all cases. Additionally, we obtain an upper estimate for the cost of null controllability. Finally, we derive a lower estimate of the cost of null controllability using a representation theorem for analytic functions of exponential type.\\ 

In this work, the initial condition $u_{0}(x)$ will be taken in the weighted Lebesgue space $L^2_\beta(0,1):=L^2((0,1);x^\beta dx)$, $\beta\in\R$, endowed with the inner product
\[\langle f,g\rangle_\beta:=\int_0^1 f(x)g(x)x^\beta dx,\]
and its corresponding norm denoted by $\|\cdot\|_{\beta}$.\\

Now, suppose that the system (\ref{problem}) admits a unique solution for each case with an initial condition in $L^{2}_{\beta}(0,1)$, a fact that will be proven in Section \ref{sec2}, although the solution lies in an interpolation space $\HH^{-s}$. The system (\ref{problem}) is said to be null controllable at time $T>0$ with controls in $L^2(0,T)$ if, for any $u_0 \in L_{\beta}^{2}(0,1)$, there exists $f \in L^2(0,T)$ such that the corresponding solution satisfies $u(\cdot,T ) \equiv 0$.

Once we establish that the system (\ref{problem}) is null controllable, we are interested in analyzing the behavior of the controllability cost. To this end, consider the set of admissible controls
\[U(T,\alpha,\beta,\mu,\ell, r, u_0)=\{f \in L^2(0, T ): u \text{ is solution of system (\ref{problem}) that satisfies } u(\cdot,T ) \equiv 0\}.\]
Then the cost of the controllability is defined as 
\[\mathcal{K}(T,\alpha,\beta,\mu, \ell, r):=\sup_{\|u_0\|_{L^{2}_{\beta}(0,1)}\leq1}\inf\left\{\|f\|_{L^2(0,T)}:f\in U(T,\alpha,\beta,\mu, \ell, r, u_0)\right\}.\]

From the system (\ref{problem}) the main result of this work is the following:
\begin{theorem}\label{Teo1}
	Let $T>0$ and $\al,\beta, \mu, \ell, \gamma \in \R$ with $0\leq\al<2$, $\mu$, $\ell$ and $\gamma$  satisfying (\ref{mucon}), (\ref{index_l}) and (\ref{gamma}) respectively, and $r\in\{0,1\}$. The next statements hold.
	\begin{enumerate}
		\item \textbf{Existence of a control.} For any $u_0\in L^2_\beta(0,1)$ there exists a control $f \in L^2(0, T )$ such that the solution $u$ to (\ref{problem}) satisfies $u(\cdot,T ) \equiv 0$.
		\item \textbf{Upper bound of the cost.} There exists a constant $c>0$ such that for every $\delta\in (0,1)$ we have that when $\gamma$ is given by (\ref{gamma}), then 
		\[\mathcal{K}(T,\alpha,\beta,\mu, \ell,r)\leq \frac{c M(T,\alpha,\nu,\delta)T^{1/2}}{\kappa_\al^{9/2-2r}\cts{\pts{\sqrt{\mu(\alpha+\beta)}+\sqrt{\mu(\alpha+\beta)-\mu}}\pts{1-\ell}+\ell}}
		\exp\pts{-\frac{T}{2}\kappa_\alpha^4 j_{\nu,1}^4}.\]
		where 
		\begin{equation}\label{Nu}
			\kappa_\al:=\frac{2-\al}{2},\quad \nu=\nu(\al,\beta,\mu):=\sqrt{\mu(\al+\beta)-\mu}/\kappa_\al,
		\end{equation}
		$j_{\nu,1}$ is the first positive zero of the Bessel function $J_\nu$ (defined in the Appendix), and
		\[
		\begin{array}{rcl}
			M(T,\alpha,\nu,\delta)& = &\d\pts{1+\frac{1}{(1-\delta)\kappa_\alpha^2 T}}\cts{\exp\pts{\frac{\sqrt{2+\sqrt{2}}}{\sqrt{2}\kappa_\alpha}}+\frac{1}{\delta^3}\exp\pts{\frac{3\sqrt{2+\sqrt{2}}}{(1-\delta)\kappa_\alpha^2 T}}}\\
			&&\\
			&&\d\times\exp\pts{-\frac{(1-\delta)^{3/2}T^{3/2}}{8\sqrt{2+\sqrt{2}}(1+T)^{1/2}}\kappa_\alpha^5 j_{\nu,1}^4}.
		\end{array}
		\]
		\item \textbf{Lower bound of the cost.} There exists a constant $c>0$ such that
		\[
		\hspace{-1cm}\begin{array}{r}
			\d\frac{c2^{\nu} \Gamma(\nu+1) \left|J_{\nu}^{\prime}\left(j_{\nu, 1}\right)\right|\exp{\left(2\pts{1-\frac{\log 5}{\pi}-\frac{\tan^{-1} 2}{\pi}}j_{\nu,2}\right)}}{\pts{{2T \kappa^{5-4r}_\alpha}}^{1/2}\cts{\pts{\sqrt{\mu(\alpha+\beta)}+\sqrt{\mu(\alpha+\beta)-\mu}}(1-\ell) + \ell}\left(j_{\nu, 1}\right)^{\nu+2-2r}}
			\exp\pts{-\pts{j_{\nu,1}^4+2j_{\nu,2}^4}\kappa_\alpha^4 T}
			\\
			\leq 
			\mathcal{K}(T,\alpha,\beta,\mu, \ell, r),
		\end{array}
		\]			
	\end{enumerate}
	where $j_{\nu,2}$ is the second positive zero of the Bessel function $J_\nu$.
\end{theorem}

The following qualitative result is derived as a consequence of Theorem \ref{Teo1}.
\begin{remark}
	Firstly, when $T\to 0^{+}$, notice that the upper bound of the control cost for the system (\ref{problem}) explodes much faster than the lower bound, since the upper bound has a behavior of the form $\frac{1}{T}\mathrm{e}^{c/T}$ with $c>0$, while the lower bound behaves as $1/\sqrt{T}$.
	
	Now, if $\alpha\to 2^{-}$, the control cost for the system blows up since $j_{\nu,1}$ satisfies (\ref{BoundFirstZero}), $j_{\nu,1} < j_{\nu,2}$ and using
	(\ref{recur}), (\ref{cerocrec}), (\ref{stirling}) and (\ref{Asymtotic}), we obtain that 
	$$
	\begin{array}{rcl}
		\d\frac{2^{\nu}\Gamma(\nu+1)\left|J_{\nu}'(j_{\nu,1})\right|}{\left(j_{\nu,1}\right)^\nu} &\sim & \d \mathrm{e}\sqrt{\frac{\nu}{\nu+1}}\left(\frac{\nu}{\nu+1}\right)^{
			\nu}\dfrac{j_{\nu,1}}{\nu+1}\\[0.5cm]
		&\geq&\d \d \mathrm{e}\sqrt{\frac{\nu}{\nu+1}}\left(\frac{\nu}{\nu+1}\right)^{
			\nu}\dfrac{\nu}{\nu+1} \rightarrow 1,\quad\text{as } \nu\to\infty.
	\end{array}
	$$
\end{remark}

This paper is organized as follows: Based on the works \cite{GaloLopez, GaloLopez2}, in Section \ref{sec2} we prove that $-\A^2$ given by (\ref{A_square}) is a negative self-adjoint operator and, therefore, a diagonalizable operator. This result enables us to consider initial data in appropriate interpolation spaces, where we introduce a notion of weak solution for each of our systems.\\

Using the moment method introduced by Fattorini \& Russell, in Section \ref{sec3} we construct a biorthogonal family involving the eigenvalues of $\A^2$. To achieve this, we employ a complex multiplier based on the decomposition in Weierstrass products of the Bessel and modified Bessel functions.\\

In Section \ref{sec4}, we prove Theorem \ref{Teo1}, establishing the controllability of the system (\ref{problem}) using the results from the previous sections. We also obtain upper estimates for the cost of null controllability for our systems. Furthermore, based on the representation Theorem (Theorem \ref{repre}), we derive a lower estimate for the cost of null controllability for each system. Finally, in Section \ref{sec5}, we discuss the case $\alpha + \beta = 1$, 
following the results of Sturm-Liouville theory as applied in \cite{GaloLopez3}, to maintain the same position as in the case $\al + \beta < 1$.

\section{Functional setting and well-posedness}\label{sec2}
In order to obtain our results, we introduce some suitable weighted spaces. For $\al,\beta\in \R$ we consider the weighted Sobolev space
\[H_{\alpha, \beta}^{1}(0,1)=\left\{u \in L_{\beta}^{2}(0,1)\cap H^1_{loc}(0,1): x^{\alpha / 2} u_{x} \in L_{\beta}^{2}(0,1)\right\}\]
with the inner product
\[\langle u, v\rangle_{\alpha, \beta}:=\int_{0}^{1} u v\, x^{\beta} \!dx+\int_{0}^{1} x^{\alpha+\beta} u_{x} v_{x} dx,\]
and its corresponding norm denoted by $\|\cdot\|_{\alpha,\beta}$.\\

Throughout this work we take $0\leq \alpha < 2$.\\ 
Then, for the case $\al+\beta < 1$, if $u\in H^{1}_{\al,\beta}(0,1)$ in \cite[Remark 2.1]{GaloLopez} was proved that the limits $\d u(0) =\lim_{x\to 0^{+}}u(x)$ and $\d u(1)=\lim_{x\to 1^{-}}u(x)$ exist and are finite.

Thus, we introduce the subspace
\begin{equation}\label{HD}
	H^{1}_{\al,\beta,0}=H^{1}_{\al,\beta,0}(0,1):=\{u\in H^{1}_{\al,\beta}(0,1)\,:\, u(0) = u(1) = 0\}.
\end{equation}

In \cite{GaloLopez} the following inequalities were proven: the generalized Hardy inequality
\begin{equation}\label{hardy}
	\mu(\al+\beta)\int_0^1\frac{u^2}{x^{2-\al-\beta}}\mathrm{d}x\leq\int_0^1x^{\al+\beta} u_x^2\mathrm{d}x,\,\,\textit{ for all } u\in H^{1}_{\al,\beta,0},
\end{equation}
where $\mu(\al+\beta)$ is defined in (\ref{mucon}) with $\delta = \al+\beta$, and the weighted Poincaré inequality
\begin{equation}\label{poincare}
	\int_{0}^{1}|u|^{2} x^{\beta} \mathrm{d}x \leq \frac{1}{(2-\alpha)(1-\alpha-\beta)} \int_{0}^{1} x^{\alpha+\beta}\left|u_{x}\right|^{2} \mathrm{d}x,\,\,\textit{ for all } u\in H^{1}_{\al,\beta,0}.
\end{equation}	
\begin{remark}
	From (\ref{hardy}) and (\ref{poincare}), it was proved in \cite[Proposition 2.5]{GaloLopez} that
	\begin{equation}\label{helpD0}
		\lim_{x\to 0^{+}}x^{\al+\beta}u_{x}(x)v(x) = 0,\,\,\lim_{x\to 1^{-}}x^{\al+\beta}u_{x}(x)v(x) = 0, \textit{ for all } u,v \in H^{1}_{\al,\beta,0}.
	\end{equation}
\end{remark}

In \cite[Propositions 2.5-2.7]{GaloLopez} was showed that $-\A:D_{0}(\A)\subset L^{2}_{\beta}(0,1) \rightarrow L^{2}_{\beta}(0,1)$ is a negative self-adjoint operator with dense domain $D_{0}(\A)$ given by
\begin{equation}\label{domainA_0}
	D_{0}(\A):=\left\{u\in H^{1}_{\al,\beta,0}\cap H^{2}_{\mathrm{loc}}(0,1)\,:\,(x^{\al}u_{x})_{x}+\beta x^{\al-1}u_{x}+\dfrac{\mu}{x^{2-\al}}u\in L^{2}_{\beta}(0,1)\right\}.
\end{equation}
Furthermore, the family
\begin{equation}\label{Phik}
	\Fi_k(x):=\frac{(2\kappa_\al)^{1/2}}{\abs{J'_\nu\pts{j_{\nu,k}}}}x^{\pts{1-\al-\beta}/2}J_\nu\pts{j_{\nu,k}x^{\kappa_\al}},\quad k\geq 1,
\end{equation}
is an orthonormal basis for $L^2_\beta(0,1)$ such that
\begin{equation}\label{lambdak}
	\mathcal{A}\Fi_k=\lambda_k \Fi_k, \quad \lambda_{k}:=\kappa_\al^2 \pts{j_{\nu,k}}^2,\quad k\geq 1,
\end{equation}
where $\nu$ is defined in (\ref{Nu}) and $\lvs{j_{\nu,k}}_{k\geq 1}$ is the sequence of positive zeros of the Bessel function $J_{\nu}$.\\

Now, we assume $\al+\beta > 1$. Then, if $u\in H^{1}_{\al,\beta}(0,1)$, from \cite[Proposition 4]{GaloLopez2} we have that the limit $u(1) = \lim_{x\to 1^{-}}u(x)$ exists and is finite. 

Therefore, we can introduce the following subspace
\begin{equation}\label{HN}
	H^{1}_{\al,\beta,N}=H^{1}_{\al,\beta,N}(0,1):=\{u\in H^{1}_{\al,\beta}(0,1)\,:\, u(1) = 0\}.
\end{equation}
Moreover, if $u\in H^{1}_{\al,\beta, N}$, then $u$ also satisfies the inequalities (\ref{hardy}) and (\ref{poincare}) (see \cite[Proposition 6 and eq. (10)]{GaloLopez2}).
\begin{remark}
	By using (\ref{hardy}) and (\ref{poincare}), in \cite[Propositions 7, 8 and 10]{GaloLopez2}  was proved that 
	\begin{equation}\label{limtsN0}
		\lim_{x\to 0^+}x^{(\al+\beta-1)/2}u(x) = \lim_{x\to 0^{+}} x^{(\al+\beta+1)/2}u_{x}(x) = 0,\,\, \text{for all } u\in H^{1}_{\al,\beta,N}.
	\end{equation}
\end{remark}
In particular, if $\A u\in H^{1}_{\al,\beta, N}$, then 
\begin{equation}\label{lim_Au}
	\lim_{x\to 0^+}x^{(\al+\beta-1)/2}\pts{\A u}(x) = \lim_{x\to 0^{+}} x^{(\al+\beta+1)/2}\pts{\A u}_{x}(x) = 0.
\end{equation}

In \cite[Propositions 10-11]{GaloLopez2} was proved that $-\A:D_{N}(\A)\subset L^{2}_{\beta}(0,1) \rightarrow L^{2}_{\beta}(0,1)$ is a negative self-adjoint operator with dense domain $D_{N}(\A)$ given by
\begin{equation}\label{domainA_N}
	D_{N}(\A):=\left\{u\in H^{1}_{\al,\beta,N}\cap H^{2}_{\mathrm{loc}}(0,1)\,:\,(x^{\al}u_{x})_{x}+\beta x^{\al-1}u_{x}+\dfrac{\mu}{x^{2-\al}}u\in L^{2}_{\beta}(0,1)\right\}
\end{equation}
and having the same orthonormal basis family given by (\ref{Phik}) for $L^2_\beta(0,1)$ that satisfies (\ref{lambdak}).

Then, from (\ref{A_square}), (\ref{HD}), (\ref{helpD0}), (\ref{domainA_0}), (\ref{lambdak}), (\ref{HN}), (\ref{limtsN0}),(\ref{lim_Au}), (\ref{domainA_N}) and the properties of $\A$ we have the following result
\begin{proposition}\label{bases}
	Let $\al,\beta\in \R$ with $0\leq \al < 2$. Then, if $\al + \beta < 1$ the operator
	\begin{equation}\label{Asquard_0}
		-\A^{2}:D_{0}(\A^2)\subset L^{2}_{\beta}(0,1) \rightarrow L^{2}_{\beta}(0,1)
	\end{equation} 
	with dense domain
	\begin{equation}\label{domainA2D}
		\begin{array}{l}
			D_{0}(\A^2):=\left\{u\in D_{0}(\A)\,:\A u\in D_{0}(\A)\right\}\\
			=\left\{u\in H^{1}_{\al,\beta,0}\cap H^4_{\mathrm{loc}}(0,1)\,:\, \cts{\A u}(0)=\cts{(x^{\al+\beta}u_{x})_{x}}(1) = 0,\,\A u,\,x^{\al/2}\pts{\A u}_{x},\,\A^2u\in L^{2}_{\beta}(0,1)\right\},
		\end{array}
	\end{equation}
	is a negative self-adjoint operator. Furthermore, the family given by (\ref{Phik})
	is an orthonormal basis for $L^2_\beta(0,1)$ such that 
	\begin{equation}\label{lambda2_k}
		\mathcal{A}^2\Fi_k=\lambda^2_k \Fi_k,\quad \lambda^2_{k} =\kappa_\al^4 \pts{j_{\nu,k}}^4,\quad k\geq 1.
	\end{equation} 
	
	In the same way, if $\al +\beta >1$, then the operator
	\begin{equation}\label{Asquard_N}
		-\A^{2}:D_{N}(\A^2)\subset L^{2}_{\beta}(0,1) \rightarrow L^{2}_{\beta}(0,1)
	\end{equation} 
	with dense domain
	
	\begin{equation}\label{domainA2N}
		\begin{aligned}
			D_N(\A^2):&=\left\{u\in D_{N}(\A)\,:\A u\in D_{N}(\A)\right\}\\
			& =\left\{u\in H^{1}_{\al,\beta,N}\cap H^{4}_{\mathrm{loc}}(0,1)\,:\, \A u\in H^{1}_{\al,\beta,N},\,\A u,\A^2u\in L^{2}_{\beta}(0,1)\right\},
		\end{aligned}
	\end{equation}
	is a negative self-adjoint operator with the same orthonormal basis family given by (\ref{Phik}) for $L^2_\beta(0,1)$ that satisfies (\ref{lambda2_k}). 
\end{proposition} 

\begin{proof}
	In view of the properties of the operator (\ref{opeA}), the results in \cite[Propositions 2.5-2.7]{GaloLopez}, \cite[Propositions 10-11]{GaloLopez2} and \cite{Hochstadt}, the result follows.
\end{proof}

Hence, in both cases, $\al+\beta < 1$ and $\al+\beta >1$, if $D(\A^2):= D_{0}(\A^2)$ or $D(\A^{2}):=D_{N}(\A^{2})$ respectively, we have that $(\A^2,D(\A^{2}))$ is the infinitesimal generator of a diagonalizable self-adjoint semigroup of contractions in $L^2_{\beta}(0,1)$. Thus, we consider interpolation spaces for the initial data. For any $s\geq 0$, we define
\[\mathcal{H}^{s}=\mathcal{H}^{s}(0,1):=D\pts{\pts{\mathcal{A}^2}^{s/2}}=\left\{u=\sum_{k=1}^\infty a_{k} \Phi_{k}:\|u\|_{\mathcal{H}^{s}}^{2}=\sum_{k=1}^\infty |a_{k}|^{2} \lambda_{k}^{2s}<\infty\right\}.\]
Also, we consider the corresponding dual spaces
\[\mathcal{H}^{-s}:=\left[\mathcal{H}^{s}(0,1)\right]^{\prime}.\]
It is well known that $\mathcal{H}^{-s}$ is the dual space of $\mathcal{H}^{s}$ with respect to the pivot space $L^2_\beta(0,1)$, i.e
\[\mathcal{H}^s\hookrightarrow \mathcal{H}^0=L^2_{\beta}(0,1)=\left(L^2_{\beta}(0,1)\right)'\hookrightarrow \mathcal{H}^{-s},\quad s>0. \]
Equivalently, $\mathcal{H}^{-s}$ is the completion of $L^2_\beta(0,1)$ with respect to the norm
\[\|u\|^2_{-s}:=\sum_{k=1}^{\infty}\lambda_k^{-2s}|\langle u,\Phi_k\rangle_\beta|^2.\]
It is well known that the linear mapping given by
\begin{equation}\label{semigroup}
	S(t)u_0=\sum_{k=1}^\infty \textrm{e}^{-\lambda^2_k t}a_k\Fi_k\quad\text{if}\quad u_0=\sum_{k=1}^\infty a_{k} \Phi_{k}\in \mathcal{H}^s
\end{equation}
defines a self-adjoint semigroup $S(t)$, $t\geq 0$, in $\mathcal{H}^s $ for all $s\in\R$.\\

For $\delta\in\R$ and a function $z:(0,1)\rightarrow \R$ we introduce the notion of $\delta$-generalized limit of $z$ at $x=0$ as follows
\[\OO_\delta (z):=\lim_{x\rightarrow 0^+} x^\delta z(x).\]

\subsection{Existence and uniqueness of solutions to system (\ref{problem})}

Now we consider a convenient definition of a weak solution for system (\ref{problem}), we multiply the equation in (\ref{problem}) by $x^\beta\varphi(\tau)=x^\beta S(\tau-t)z^{\tau}$, integrate by parts (formally), hence
\begin{equation}\label{partes}
	\begin{array}{rcl}
		\d\int_{0}^{1} \left.\left[ x^{\beta} u\varphi\right] \right|_{t=0}^{t=\tau}dx & = &\d-\int_{0}^{\tau}\left.\left[x^{\alpha+\beta} \varphi \left(\A u\right)_{x}\right]\right|_{x=0}^{x=1}dt + \int_{0}^{\tau} \left.\left[x^{\alpha+\beta} \varphi_{x} \A u\right]\right|_{x=0}^{x=1}dt\\[0.3cm]
		& &\d - \int_{0}^{\tau} \left.\left[x^{\alpha+\beta} u\A\varphi\right]\right|_{x=0}^{x=1}dt +\int_{0}^{\tau} \left.\left[x^{\alpha+\beta} \left(\A \varphi\right)_{x}u\right]\right|_{x=0}^{x=1}dt.
	\end{array}
\end{equation}
Now, we will consider the different boundary conditions for each case $\alpha + \beta <1$ and $\alpha + \beta >1$ respectively to derive the following definition of weak solution.

\begin{definition}
	Let $T>0$ and $\al,\beta,\mu, r\in \R$ with $0\leq\al<2$, $\mu<\mu(\alpha+\beta)$ and $r\in\{0,1\}$. Let $f \in L^2(0,T)$ and $u_0\in \HH^{-s}$ for some $s > 0$. A weak solution of (\ref{problem}) is a function $u \in C^0([0,T];\HH^{-s})$ such that for every $\tau \in (0,T]$ and for
	every $z^\tau \in \HH^s$ we have
	\begin{equation}\label{weaksol}
		\left\langle u(\tau), z^{\tau}\right\rangle_{\mathcal{H}^{-s}, \mathcal{H}^{s}}=\left\langle u_{0}, S(\tau)z^\tau\right\rangle_{\mathcal{H}^{-s}, \mathcal{H}^{s}} + (-1)^{1-\ell}\int_{0}^{\tau} f(t) \mathcal{O}_{\al+\beta+\gamma-\ell}\pts{\partial^{1-\ell}_{x}\cts{\A^{1-r} S(\tau-t)z^{\tau}}}\mathrm{d} t,
	\end{equation}
	where $\ell=\ell(\al,\beta)$ and $\gamma=\gamma(\alpha,\beta,\mu)$ are given by (\ref{index_l}) and (\ref{gamma}) respectively, and $A^{0}$ is the identity operator.
\end{definition}
The next result show the existence of weak solutions for the system (\ref{problem}) under suitable conditions on the parameters $\alpha,\beta,\mu,\ell, r, \gamma$ and $s$.
\begin{proposition}\label{continuity}
	Let $T>0$ and $\al,\beta,\mu, \ell,r,\gamma\in \R$ given as in the last definition. Let $f \in L^2(0,T)$ and $u_0\in \HH^{-s}$ such that $s>(\nu+1)/2-r$, where $\nu$ is given in (\ref{Nu}). Then, formula (\ref{weaksol}) defines for each $\tau \in [0, T ]$ a unique element $u(\tau) \in \HH^{-s}$ that can be written as
	\[
	u(\tau)=S(\tau) u_{0}+(-1)^{1-\ell}B_{\ell r}(\tau) f, \quad \tau \in(0, T],\]
	where $B_{\ell r}(\tau)$ is the strongly continuous family of bounded operators $B_{\ell r}(\tau): L^{2}(0,T) \rightarrow \mathcal{H}^{-s}$ given by
	\[\left\langle B_{\ell r}(\tau) f, z^{\tau}\right\rangle_{\mathcal{H}^{-s}, \mathcal{H}^{s}}=\int_{0}^{\tau} f(t) \mathcal{O}_{\al+\beta+\gamma-\ell}\pts{\partial^{1-\ell}_{x}\cts{\A^{1-r} S(\tau-t)z^{\tau}}}\mathrm{d} t, \quad \text{for all  }z^{\tau} \in \mathcal{H}^{s} .\]
	Furthermore, the unique weak solution $u$ on $[0, T]$ to (\ref{problem}) (in the sense of (\ref{weaksol})) belongs to ${C}^{0}\left([0, T] ; \mathcal{H}^{-s}\right)$ and fulfills
	\[
	\|u\|_{L^{\infty}\left([0, T] ; \mathcal{H}^{-s}\right)} \leq C\left(\left\|u_{0}\right\|_{\mathcal{H}^{-s}}+\|f\|_{L^{2}(0, T)}\right).
	\]
\end{proposition}
\begin{proof}
	Set $\tau>0$, $r\in\{0,1\}$ and $\ell$ given by (\ref{index_l}). Let $u(\tau)\in H^{-s}$ be determined by the condition (\ref{weaksol}), hence
	$$
	(-1)^{1-\ell}u(\tau)+(-1)^{\ell}S(\tau) u_{0}=\zeta_{\ell r}(\tau)f,
	$$
	where
	$$
	\left\langle\zeta_{\ell r}(\tau)f, z^{\tau}\right\rangle_{\mathcal{H}^{-s}, \mathcal{H}^{s}}=\int_{0}^{\tau} f(t) \mathcal{O}_{\al+\beta+\gamma-\ell}\pts{\partial^{1-\ell}_{x}\cts{\A^{1-r} S(\tau-t)z^{\tau}}}\mathrm{d} t, \quad \text{for all  } z^{\tau} \in \mathcal{H}^{s}.
	$$
	We claim that $\zeta_{\ell r}(\tau)$ is a bounded operator from $L^{2}(0, T)$ into $\mathcal{H}^{-s}$: Let $z^{\tau} \in \mathcal{H}^{s}$ given by 
	\begin{equation}\label{finaldata}
		z^{\tau}=\sum_{k=1}^\infty a_{k} \Phi_{k},
	\end{equation}
	therefore, if we write $d^{1-j}\Fi_k/dx^{1-j}=\Fi^{(1-j)}_k$, from (\ref{semigroup}) we have that
	$$
	\partial^{1-\ell}_{x}\cts{\A^{1-r} S(\tau-t)z^{\tau}}=\sum_{k=1}^\infty \lambda^{1-r}_{k}\mathrm{e}^{\lambda^{2}_{k}(t-\tau)} a_{k} \Phi^{(1-\ell)}_{k}, \quad \text{for all } t \in[0, \tau],
	$$
	
	By using Lemma \ref{reduce} and (\ref{limAux}) we obtain that there exists a constant $C=C(\al,\beta,\mu)>0$ such that
	\begin{equation}\label{cotaODerPhik}
		|\mathcal{O}_{\alpha+\beta+\gamma-\ell}\left(\Fi^{(1-\ell)}_k\right)|\leq C |j_{\nu,k}|^{\nu+1/2},\quad k\geq 1.
	\end{equation}
	Hence (\ref{below}) implies that there exists $C=C(\al,\beta,\mu)>0$ such that
	$$
	\begin{array}{rl}
		\hspace{0.5cm}\d\pts{\int_{0}^{\tau}\left|\mathcal{O}_{\alpha+\beta+\gamma-\ell}\left(\partial^{1-\ell}_{x}\cts{\A^{1-r} S(\tau-t)z^{\tau}}\right)\right|^{2} \mathrm{~d} t}^{1/2}&\\ 
		&\hspace{-2.5cm}
		\d\leq  \sum_{k=1}^\infty 
		|a_k||\lambda^{1-r}_{k}||\OO_{\al+\beta+\gamma-\ell}(\Fi^{(1-r)}_k)| \pts{\int_0^\tau \mathrm{e}^{2\lambda^{2}_{k}(t-\tau)} \mathrm{~d} t}^{1/2} \\
		&\hspace{-2.5cm}
		\d\leq C\left\|z^{\tau}\right\|_{\mathcal{H}^{s}}\pts{\sum_{k=1}^\infty |\lambda_k|^{\nu+5/2-2s-2r}\int_0^\tau \mathrm{e}^{2\lambda^{2}_{k}(t-\tau)} \mathrm{~d} t}^{1/2}\\
		&\hspace{-2.5cm}
		\d = C\left\|z^{\tau}\right\|_{\mathcal{H}^{s}}\pts{\sum_{k=1}^\infty |\lambda_k|^{\nu+1/2-2s-2r}\pts{1-\mathrm{e}^{-2\lambda^{2}_k \tau}}}^{1/2}\\
		&\hspace{-2.5cm}
		\leq \d C\left\|z^{\tau}\right\|_{\mathcal{H}^{s}}\pts{\sum_{k=1}^\infty\frac{1}{k^{2(2s-\nu-1/2+2r)}}}^{1/2}=C\left\|z^{\tau}\right\|_{\mathcal{H}^{s}}.
	\end{array}
	$$
	Therefore $\|\zeta_{\ell r}(\tau) f\|_{\mathcal{H}^{-s}}\leq C\|f\|_{L^2(0,T)}$ for all $f\in L^2(0,T)$, $\tau\in (0,T]$.\\
	
	Finally, we fix $f\in L^2(0,T)$ and show that the mapping $\tau\mapsto \zeta_{\ell r}(\tau) f$ is right-continuous on $[0,T)$. Let $h>0$ small enough and $z\in \mathcal{H}^s$ given as in (\ref{finaldata}). Thus, proceeding as in the last inequalities, we have
	$$
	\begin{array}{rl}
		\d|\left\langle\zeta_{\ell r}(\tau+h)f-\zeta_{\ell r}(\tau)f, z\right\rangle_{\mathcal{H}^{-s}, \mathcal{H}^{s}}|&
		\\[0.3cm]
		&\hspace{-2.9cm}
		\d\leq\int_{0}^{\tau} |f(t)| |\mathcal{O}_{\alpha+\beta+\gamma-\ell}\left(\partial^{1-\ell}_{x}\cts{\A^{1-r}\pts{S(\tau+h-t) -S(\tau-t)}z}\right)| \mathrm{d} t \\
		& \hspace{0.7cm}
		\d + \int_{\tau}^{\tau+h} |f(t)| |\mathcal{O}_{\alpha+\beta+\gamma-\ell}\left(\partial_{x}^{1-\ell}\cts{\A^{1-r} S(\tau+h-t) z}\right)| \mathrm{d} t\\[0.3cm]
		&\hspace{-2.9cm}
		\d\leq C\left\|z\right\|_{\mathcal{H}^{s}}\|f\|_{L^2(0,T)}\cts{\pts{\sum_{k=1}^\infty\frac{I(\tau,k,h)}{k^{2(2s-\nu-1/2+2r)}}}^{1/2}+\pts{\sum_{k=1}^\infty\frac{1-\mathrm{e}^{-2\lambda^{2}_k h}}{k^{2(2s-\nu-1/2+2r)}}}^{1/2}},
	\end{array}
	$$
	where 
	\begin{equation}\label{limI(t,k,h)}
		\begin{array}{rcl}
			I(\tau,k,h)&=&\d\lambda^{2}_k\int_0^\tau\pts{\mathrm{e}^{\lambda^{2}_k(t-\tau-h)}-\mathrm{e}^{\lambda^2_k(t-\tau)}}^2\mathrm{~d} t\\
			&=&\d \frac{1}{2}(1-\mathrm{e}^{-\lambda^{2}_k h})^2(1-\mathrm{e}^{-2\lambda^{2}_k \tau})\rightarrow 0\quad\text{as}\quad h\rightarrow 0^+.
		\end{array}
	\end{equation}
	Since $0\leq I(\tau,k,h)\leq 1/2$ uniformly for $\tau, h>0$, $k\geq 1$, the result follows by the dominated convergence theorem.
\end{proof}

\begin{remark}
	In the following sections, we will consider initial conditions in $L^2_\beta(0,1)$. Notice that $L^2_\beta(0,1)\subset H^{-(\nu+1+\delta-2r)/2}$ for all $\delta> \min\{2r-\nu-1, 0\}$, and we can apply Proposition \ref{continuity} with $s=(\nu+1+\delta-2r)/2$, $\delta>0$, then the corresponding solutions will be in $C^0\pts{[0, T ]; H^{-(\nu+1+\delta-2r)/2}}$.
\end{remark}

\section{The biorthogonal family and its multiplier}\label{sec3}
In this section we obtain some tools to prove the null controllability of each system based on the ideas given by \cite{GaloLopez,GaloLopez2,GaloLopez3}, where it was applied the method of moments introduced by Fattorini \& Russell in \cite{Fattorini}, which in our case involves the construction of a family $\d \lvs{\psi_k}_{k\geq 1}\subset L^2(0,T)$ biorthogonal to the family of exponential functions $\lvs{\mathrm{e}^{-\lambda^2_{k}(T-t)}}_{k\geq 1}$ on $[0, T]$, i.e that it satisfies
$$
\int_{0}^{T}\psi_k(t) \mathrm{e}^{-\lambda^2_{l}(T-t)} dt = \delta_{kl},\quad\text{for all}\quad k,l\geq 1,\,\text{ and } \d\delta_{k\ell}:=\left\{\begin{array}{rl} 1, & \text{ if } k =\ell,\\ 0, & \text{ if } k \neq \ell.\end{array}\right.
$$

As a consequence, we will get an upper bound for the cost of the null controllability of our systems.\\

Assume that each $F_k$, $k\geq 1$, is an entire function of exponential type $T/2$ such that $F_k\in L^2(\R)$, and
\begin{equation}\label{krone}
	F_{k}(i\lambda^2_{l})=\delta_{kl},\quad \text{for all}\quad k,l\geq 1.
\end{equation}
The $L^2$-version of the Paley-Wiener theorem implies that there exists $\eta_k\in L^2(\R)$ with support in $[-T/2,T/2]$ such that $F_k(z)$ is the analytic extension of the Fourier transform of $\eta_k$. Then we have that 
\begin{equation}\label{psieta}
	\psi_k(t):=\mathrm{e}^{\lambda^2_k T/2}\eta_k(t-T/2),\quad t\in[0,T],\,k\geq1,
\end{equation}
is the family we are looking for:
\[\delta_{kl}=\mathrm{e}^{(\lambda^2_k-\lambda^2_l)T/2}F_{k}(i\lambda^2_{l})=\mathrm{e}^{(\lambda^2_k-\lambda^2_l)T/2}\int_{-\frac{T}{2}}^{\frac{T}{2}} \eta_{k}(t) \mathrm{e}^{\lambda^2_{l}t} \mathrm{d}t=\int_{0}^{T}\psi_k(t) \mathrm{e}^{-\lambda^2_{l}(T-t)} \mathrm{d}t\quad\text{for all}\quad k,l\geq 1.\]

Now, we proceed to construct the family $F_k$, $k\geq 1$, satisfying the aforementioned properties. First, consider the Weierstrass infinite product
\begin{equation}
	\Lambda(z):=\prod_{k=1}^{\infty}\pts{1+\dfrac{iz}{\pts{\kappa_\al j_{\nu, k}}^4}}.
\end{equation}
By using (\ref{asint}) we have that $j_{\nu, k}=O(k)$ for $k$ large, thus the infinite product is well defined and converges absolutely in $\C$. Hence, $\Lambda(z)$ is an entire function with simple zeros at $i\pts{\kappa_\al j_{\nu, k}}^4=i\lambda^2_k$, $k\geq 1$.\\

Using \cite[Chap. XV, p. 498, eq. (3)]{Watson}, we can write
\begin{equation}\label{fLambda}
	\begin{array}{rcl}
		\Lambda(z) & = &\cts{\Gamma\pts{\nu+1}}^2\pts{\dfrac{4\kappa^2_\al}{i\sqrt{-i\,z}}}^\nu J_{\nu}\pts{\dfrac{\sqrt[4]{-iz}}{\kappa_\al}}J_{\nu}\pts{\dfrac{i\sqrt[4]{-iz}}{\kappa_\al}}\\
		&=& \cts{\Gamma\pts{\nu+1}}^{2}\pts{\dfrac{4\kappa^2_\al}{\sqrt{-i\,z}}}^\nu J_{\nu}\pts{\dfrac{\sqrt[4]{-iz}}{\kappa_\al}}I_{\nu}\pts{\dfrac{\sqrt[4]{-iz}}{\kappa_\al}},
	\end{array}
\end{equation}
where $I_{\nu}(z) = i^{-\nu}J_{\nu}(iz)$ is the modified Bessel function (see \cite[Chap. III, p.77]{Watson}).
In \cite[pag. 12]{GaloLopez} was proved that
\[|J_\nu(z)|\leq \frac{|z|^\nu \mathrm{e}^{|\Im(z)|}}{2^\nu\Gamma\left(\nu+1\right)},\quad z\in\C.\]
Therefore,
\[|\Lambda(z)|\leq \exp\pts{\frac{\left|\Im\pts{\sqrt[4]{-iz}}\right|+\left|\Im\pts{i\sqrt[4]{-iz}}\right|}{\kappa_\alpha}},\quad z\in\C.\]
In particular,
\begin{equation}\label{besst}
	|\Lambda(z)|\leq\exp\pts{\frac{\sqrt{2}\,|z|^{1/4}}{\kappa_\al}},\quad z\in\C,\quad |\Lambda(x)|\leq\exp\pts{\frac{\sqrt{2+\sqrt{2}}|x|^{1/4}}{\sqrt{2}\kappa_\al}},\quad x\in\R.
\end{equation}

It follows that 
\begin{equation}\label{PsiFunction}
	\Psi_{k}(z):=\dfrac{\Lambda(z)}{\Lambda'(i\lambda^2_{k})(z-i\lambda^2_{k})},\quad k\geq 1,
\end{equation}
is a family of entire functions that satisfy (\ref{krone}). Since $\Psi_{k}(x)$ is not in $L^2(\R)$, we need to fix this by using a suitable ``complex multiplier", thus we follow the approach introduced in \cite{Tucsnak}.\\

For $\theta, a>0$, we define
$$
\sigma_{\theta}(t):=\exp\pts{-\frac{\theta}{1-t^2}},\quad t\in(-1,1),
$$
and extended by $0$ outside of $(-1, 1)$. Clearly $\sigma_{\theta}$ is analytic on $(-1,1)$. Set $C_{\theta}^{-1}:=\int_{-1}^{1}\sigma_{\theta}(t)\mathrm{d}t$ and define
\begin{equation}\label{Hfunction}
	H_{a,\theta}(z)=C_{\theta}\int_{-1}^{1}\sigma_{\theta}(t)\exp\pts{-ia tz}\mathrm{d}t.
\end{equation}
$H_{a,\theta}(z)$ is an entire function, and the next result provides additional properties of $H_{a,\theta}(z)$.\\
\begin{lemma}
	The function $H_{a,\theta}$ fulfills the following inequalities
	\begin{eqnarray}
		H_{a,\theta}(ix)&\geq &\frac{\exp\pts{a|x|/\pts{2\sqrt{\theta+1}}}}{11\sqrt{\theta+1}},\quad x\in\R,\label{Hcot1}\\
		|H_{a,\theta}(z)|     &\leq & \exp\pts{a|\Im(z)|},\quad z\in\C,\label{Hcot2}\\
		|H_{a,\theta}(x)|     &\leq & \chi_{|x|\leq 1}(x)+c\sqrt{\theta+1}\sqrt{a\theta\abs{x}}\exp\pts{3\theta/4-\sqrt{a\theta\abs{x}}}\chi_{|x|> 1}(x),\quad x\in\R,\label{Hcot3}
	\end{eqnarray}
	where $c>0$ does not depend on $a$ and $\theta$.
\end{lemma}
We refer to \cite[pp. 85--86]{Tucsnak} for the details.\\

Finally, for $k\geq 1$ consider the entire function $F_{k}$ given as
\begin{equation}\label{Ffunction}
	F_{k}(z):=\Psi_{k}(z)\dfrac{H_{a,\theta}(z)}{H_{a,\theta}(i\lambda^2_{k})},\quad z\in\C.
\end{equation}

For $\delta\in(0,1)$ we set
\begin{equation}\label{aConst}
	a:=\frac{T(1-\delta)}{2}>0,\quad \text{and}\quad \theta:=\dfrac{\pts{2+\sqrt{2}}(1+\delta)^2}{\kappa_\al^2 T\pts{1-\delta}}>0.
\end{equation}

The computations in the next result justify the choice of these parameters.
\begin{lemma}
	For each $k\geq 1$ the function $F_{k}(z)$ satisfies the following properties:\\
	i) $F_{k}$ is of exponential type $T/2$,\\
	ii) $F_{k}\in L^1(\R)\cap L^2(\R)$,\\
	iii) $F_k$ satisfies (\ref{krone}).\\
	iv) Furthermore, there exists a constant $c>0$, independent of $T,\alpha$ and $\delta$, such that
	\begin{equation}\label{Fbound}
		\left\|F_{k}\right\|_{L^{1}(\R)} \leq \frac{C(T, \alpha,\delta)}{\lambda^2_{k}\left|\Lambda^{\prime}\left(i \lambda^2_{k}\right)\right|}  \exp\pts{-\frac{a\lambda^2_k}{2\sqrt{\theta+1}}},		
	\end{equation} 
	where
	\begin{equation}\label{upper}
		C(T, \alpha,\delta)=c\sqrt{\theta+1}\cts{\exp\pts{{\frac{\sqrt{2+\sqrt{2}}}{\sqrt{2}\kappa_\alpha}}}+\sqrt{\theta+1}\frac{\kappa_\alpha^2}{\delta^3}\exp\pts{\frac{3 \theta}{4}}}. 
	\end{equation}
\end{lemma}
\begin{proof}
	By using (\ref{besst}), (\ref{Hcot2}), (\ref{Ffunction}) and (\ref{aConst}) we get that $F_{k}$ is of exponential type $T/2$. Moreover, by (\ref{PsiFunction}) and (\ref{Ffunction}), $F_{k}$ fulfills (\ref{krone}).\\
	
	Now we use (\ref{besst}), (\ref{Hcot1}), (\ref{Hcot3}), (\ref{Ffunction}), and (\ref{aConst}) to get
	$$
	\begin{array}{rl}
		\d\left|F_{k}(x)\right| & \d \leq c\exp\pts{-\frac{a\lambda^2_k}{2\sqrt{\theta+1}}}\frac{\sqrt{\theta+1}}{|\Lambda^{\prime}\left(i \lambda^2_{k}\right)||x^2+ \lambda_{k}^4|^{1/2}} |H_{a,\theta}(x)|\exp\pts{\frac{\sqrt{2+\sqrt{2}}|x|^{1/4}}{\sqrt{2}\kappa_\al}}
		\\[0.3cm]
		& \d\leq c \exp\pts{-\frac{a\lambda^2_k}{2\sqrt{\theta+1}}}\frac{\sqrt{\theta+1}}{\lambda^2_{k}|\Lambda^{\prime}\left(i \lambda^2_{k}\right)|}\\[0.3cm]
		&\hspace{1cm}\d\times\cts{\mathrm{e}^{\frac{\sqrt{2+\sqrt{2}}}{\sqrt{2}\kappa_\alpha}}\chi_{|x|\leq 1}(x)+\sqrt{\theta+1}\sqrt{a\theta\abs{x}}\exp\pts{\frac{3\theta}{4}- \frac{\delta\sqrt{2+\sqrt{2}}|x|^{1/2}}{\sqrt{2}\kappa_\alpha}}\chi_{|x|> 1}(x)}.
	\end{array}
	$$	
	Since the function on the right-hand side is rapidly decreasing in $\R$, we have $F_{k}\in L^1(\R)\cap L^2(\R)$. Finally, the change of variable $y=\delta\sqrt{2+\sqrt{2}}\,|x|^{1/2}/\pts{\sqrt{2}\kappa_\al}$ implies (\ref{Fbound}).
\end{proof}

Since $\eta_k, F_k\in L^1(\R)$, the inverse Fourier theorem yields 
\[\eta_k(t)=\frac{1}{2\pi}\int_{\R}\mathrm{e}^{it\tau}F_k(\tau)\mathrm{d}\tau,\quad t\in\R, k\geq 1,\]
hence (\ref{psieta}) implies that $\psi_k\in C([0,T])$ and from (\ref{Fbound}) we have
\begin{equation}\label{psiL1}
	\|\psi_k\|_{\infty}\leq \frac{C(T, \alpha,\delta)}{\lambda^2_{k}\left|\Lambda^{\prime}\left(i \lambda^2_{k}\right)\right|}  \exp\pts{\frac{\lambda^2_k T}{2}-\frac{a\lambda^2_k}{2\sqrt{\theta+1}}},\quad k\geq 1.
\end{equation}

\section{Proof of Theorem \ref{Teo1}: Control at the left endpoint}\label{sec4}
\subsection{Upper estimate of the cost of the null controllability}
We are ready to prove the null controllability of our first system (\ref{problem}). Let $u_{0}\in L^{2}_\beta(0,1)$, then consider its Fourier series with respect to the orthonormal basis $\{\Phi_{k}\}_{k\geq 1}$,
\begin{equation}\label{uoSerie}
	u_{0}(x)=\sum_{k=1}^{\infty} a_{k} \Phi_{k}(x).
\end{equation}
We set
\begin{equation}\label{f1D0serie}
	f(t):=\sum_{k=1}^{\infty}\frac{(-1)^{\ell}a_{k} \mathrm{e}^{-\lambda^2_{k} T}}{\lambda^{1-r}_{k}\mathcal{O}_{\alpha+\beta+\gamma-\ell}\left(\Phi^{(1-\ell)}_{k}\right)} \psi_{k}(t),
\end{equation}
where $r\in\{0,1\}$, $\ell$ and $\gamma$ are given by (\ref{index_l}) and (\ref{gamma}), respectively. Since $\{\psi_k\}$ is biorthogonal to $\{\mathrm{e}^{-\lambda^2_k(T-t)}\}$, we have
$$
\begin{array}{rcl}
	\d\int_{0}^{T} f(t) \lambda^{1-r}_{k}\mathcal{O}_{\alpha+\beta+\gamma-\ell}\left(\Phi^{(1-\ell)}_{k}\right) \mathrm{e}^{-\lambda^2_{k}(T-t)} \mathrm{d} t &=&\d (-1)^{\ell}a_{k} \mathrm{e}^{-\lambda^2_{k} T}\\ 
	&=& \d(-1)^{\ell}\left\langle u_{0}, \mathrm{e}^{-\lambda^2_{k} T}\Phi_{k}\right\rangle_\beta\\ 
	&=& \d (-1)^{\ell}\left\langle u_{0}, \mathrm{e}^{-\lambda^2_{k} T}\Phi_{k}\right\rangle_{\mathcal{H}^{-s}, \mathcal{H}^{s}}.
\end{array}
$$
Let $u\in C([0,T];H^{-s})$ that satisfies (\ref{weaksol}) for all $\tau\in (0,T]$, $z^\tau\in H^s$. In particular, for $\tau=T$ we take $z^T=\Phi_k$, $k\geq 1$, then the last equality implies that
$$
\left\langle u(\cdot, T), \Phi_{k}\right\rangle_{\mathcal{H}^{-s}, \mathcal{H}^{s}}=0\quad \text{for all}\quad k \geq 1,
$$
hence $u(\cdot, T)=0$.\\

It just remains to estimate the norm of the control $f$. From (\ref{psiL1}) and (\ref{f1D0serie})  we get
\begin{equation}\label{finty}
	\|f\|_{\infty} \leq  C(T, \alpha,\delta)\sum_{k=1}^{\infty} \frac{\left|a_{k}\right|}{\left|\mathcal{O}_{\alpha+\beta+\gamma-\ell}\left(\Phi^{(1-\ell)}_{k}\right)\right|} \frac{1}{\lambda^{3-r}_{k}\left|\Lambda^{\prime}\left(i \lambda^2_{k}\right)\right|} \exp\pts{-\frac{T\lambda^2_k}{2}-\frac{a\lambda^2_k}{2\sqrt{\theta+1}}}.
\end{equation}

From (\ref{fLambda}), we obtain
\begin{equation}\label{fLambdaDer}
	\left|\Lambda^{\prime}\left(i \lambda^2_{k}\right)\right|=\cts{\Gamma(\nu+1)}^2\frac{4^{\nu-1}}{\kappa^2_\al|j_{\nu, k}|^{2\nu+3}} |J_{\nu}^{\prime}\left(j_{\nu, k}\right)||I_{\nu}\left(j_{\nu, k}\right)|, \quad k\geq 1,
\end{equation}
and by using (\ref{lambdak}) and (\ref{limAux}) we get
\begin{equation}\label{estimate}
	\d\left|\mathcal{O}_{\alpha+\beta+\gamma-\ell}\left(\Phi^{(1-\ell)}_{k}\right)\lambda^{3-r}_{k}\Lambda^{\prime}\left(i \lambda^2_{k}\right)\right|= 2^{\nu-3/2}\kappa^{9/2-2r}_\al\Gamma(\nu+1)\frac{\pts{\sqrt{\mu(\alpha+\beta)}+\kappa_\al\nu}(1-\ell)+\ell}{|j_{\nu,k}|^{\nu-3+2r}}| I_{v}\pts{j_{\nu, k}}|.
\end{equation}
From (\ref{finty}), combining (\ref{estimate}), (\ref{lowerBM}) and using that $\lambda_k\geq \lambda_1$, it follows that 
\begin{equation*}
	\|f\|_{\infty} \leq \frac{C(T, \alpha,\delta)}{\kappa_\al^{9/2-2r}\cts{\pts{\sqrt{\mu(\alpha+\beta)}+\kappa_\al\nu}(1-\ell)+\ell}} \exp\pts{-\frac{T\lambda^2_1}{2}-\frac{a\lambda^2_1}{2\sqrt{\theta+1}}}\sum_{k=1}^{\infty} \frac{|a_{k}| }{j^{3-2r}_{\nu, k}}.
\end{equation*}
By using the Cauchy-Schwarz inequality, the fact that $j_{\nu,k}\geq (k-1/4)\pi$ (by (\ref{below}))  and (\ref{uoSerie}), we obtain that
\begin{eqnarray*}
	\|f\|_{\infty} &\leq& \frac{C(T, \alpha,\delta)}{\kappa_\al^{9/2-2r}\cts{\pts{\sqrt{\mu(\alpha+\beta)}+\kappa_\al\nu}(1-\ell)+\ell}} \exp\pts{-\frac{T\lambda^2_1}{2}-\frac{a\lambda^2_1}{2\sqrt{\theta+1}}}\left(\sum_{k=1}^{\infty} |a_{k}|^{2}\right)^{\frac{1}{2}}\\
	&=&\frac{C(T, \alpha,\delta)}{\kappa_\al^{9/2-2r}\cts{\pts{\sqrt{\mu(\alpha+\beta)}+\kappa_\al\nu}(1-\ell)+\ell}} \exp\pts{-\frac{T\lambda^2_1}{2}-\frac{a\lambda^2_1}{2\sqrt{\theta+1}}}\left\|u_{0}\right\|_{\beta}.
\end{eqnarray*}
Using the expression of $a,\theta$ given in (\ref{aConst}) and the facts $\theta>0$, $\delta\in(0,1)$, and $0<\kappa_\alpha\leq 1$, we get that
\[\theta\leq \frac{4\sqrt{2+\sqrt{2}}}{(1-\delta)\kappa_\alpha^2 T },\quad \sqrt{\theta+1}\leq \frac{2\sqrt{2+\sqrt{2}}(1+T)^{1/2}}{(1-\delta)^{1/2}\kappa_\alpha T^{1/2}},\quad \sqrt{\theta+1}\leq \theta+1,\]
therefore
\[
\frac{a}{\sqrt{\theta+1}}\geq \frac{\kappa_\alpha(1-\delta)^{3/2}T^{3/2}}{4\sqrt{2+\sqrt{2}}(1+T)^{1/2}}, 
\]
\[C(T, \alpha,\delta)\leq c\pts{1+\frac{1}{(1-\delta)\kappa_\alpha^2 T}}\cts{\exp\pts{\frac{\sqrt{2+\sqrt{2}}}{\sqrt{2}\kappa_\alpha}}+\frac{1}{\delta^3}\exp\pts{\frac{3\sqrt{2+\sqrt{2}}}{(1-\delta)\kappa_\alpha^2 T}}},\]
and by using the definition of $\lambda_{1}$ the result follows. 

\subsection{Lower estimate of the cost of the null controllability}
In this section we get a lower estimate of the cost $\mathcal{K}=\mathcal{K}(T,\alpha,\beta,\mu, \ell, r)$. We set
\begin{equation}\label{u0first} 
	u_0(x):= \frac{\abs{J'_\nu\pts{j_{\nu,1}}}}{(2\kappa_\al)^{1/2}}\Fi_1(x),\,x\in(0,1),
	\quad \text{hence}\quad\|u_0\|^2_\beta=\frac{\abs{J'_\nu\pts{j_{\nu,1}}}^2}{2\kappa_\al}.
\end{equation}
For $\varepsilon>0$, there exists $f\in U(T, \al,\beta,\mu, \ell, r, u_0)$ such that{\tiny }
\begin{equation}\label{inicost}
	u(\cdot,T)\equiv 0,\quad \text{and}\quad \|f\|_{L^2(0,T)}\leq (\mathcal{K}+\varepsilon)\|u_0\|_\beta.
\end{equation}
Then, in (\ref{weaksol}) we set $\tau=T$ and  take $z^\tau=\Fi_k$, $k\geq 1$,  to obtain
\begin{eqnarray*}
	\mathrm{e}^{-\lambda^2_k T}\left\langle u_{0},\Fi_k\right\rangle_\beta=\left\langle u_{0}, S(T)\Fi_k\right\rangle_{\mathcal{H}^{-s}, \mathcal{H}^{s}}& = &
	(-1)^{\ell}\int_{0}^{T} f(t) \mathcal{O}_{\al+\beta+\gamma-\ell}\left(\partial^{1-\ell}_{x}\cts{\A ^{1-r}S(T-t)\Fi_k}\right) \mathrm{d} t\\
	& = &
	(-1)^{\ell}\lambda^{1-r}_{k}\mathrm{e}^{-\lambda^2_k T}\mathcal{O}_{\al+\beta+\gamma-\ell}\left(\Fi^{(1-\ell)}_k\right)\int_{0}^{T} f(t) \mathrm{e}^{\lambda^2_k t} \mathrm{d} t,
\end{eqnarray*}
from (\ref{u0first}) and (\ref{limAux}) it follows that
\begin{equation}\label{orto}
	\int_{0}^{T} f(t) \mathrm{e}^{\lambda^2_k t} \mathrm{d} t= \frac{(-1)^{\ell}2^{\nu}\Gamma(\nu+1)\abs{J'_{\nu}\pts{j_{\nu, 1}}}^2}{2\kappa^{3-2r}_\al\cts{\pts{\sqrt{\mu(\alpha+\beta)}+\kappa_\al \nu}(1-\ell)+\ell}\pts{j_{\nu,1}}^{\nu+2-2r}}\delta_{1,k},\quad k\geq 1.
\end{equation}

Now consider the function $v: \mathbb{C} \rightarrow \mathbb{C}$ given by
\begin{equation}\label{vanalitica}
	v(s):=\int_{-T / 2}^{T / 2} f\left(t+\frac{T}{2}\right) \mathrm{e}^{-i s t} \mathrm{~d} t, \quad s \in \mathbb{C} .
\end{equation}
Fubini and Morera's theorems imply that $v(s)$ is an entire function. Moreover, (\ref{orto}) implies that
\[v(i\lambda^2_k)=0\,\,\text{for all }k\geq 2,\,\, \text{and}\,\, v(i\lambda^2_1)= \frac{(-1)^{\ell}2^{\nu}\Gamma(\nu+1)\abs{J'_{\nu}\pts{j_{\nu, 1}}}^2}{2\kappa^{3-2r}_\al\cts{\pts{\sqrt{\mu(\alpha+\beta)}+\kappa_\al \nu}(1-\ell)+\ell}\pts{j_{\nu,1}}^{\nu+2-2r}}\mathrm{e}^{-\lambda^{2}_1 T/2}.\]
We also have that

\begin{equation}\label{uve}
	|v(s)| \leq \mathrm{e}^{T|\Im(s)|/2}\int_{0}^{T}|f(t)| \mathrm{d} t \leq (\mathcal{K}+\varepsilon) T^{1/2}\mathrm{e}^{T|\Im(s)|/2} \left\|u_{0}\right\|_{\beta}.
\end{equation}
Consider the entire function $F$ given by
\begin{equation}\label{entire}
	F(s):=v\left(s-i \delta\right), \quad s \in \mathbb{C},
\end{equation}
for some $\delta>0$ that it will be chosen later on. Let $a_{k}:=i\left(\lambda^2_k+\delta\right),\quad k\geq 1$, then we have
\begin{equation}\label{ena}
	F\left(a_{k}\right)=0, \quad k\geq 2,\quad\text{and}\quad
	F\left(a_{1}\right)= \frac{(-1)^{\ell}2^{\nu}\Gamma(\nu+1)\abs{J'_{\nu}\pts{j_{\nu, 1}}}^2}{2\kappa^{3-2r}_\al\cts{\pts{\sqrt{\mu(\alpha+\beta)}+\kappa_\al \nu}(1-\ell)+\ell}\pts{j_{\nu,1}}^{\nu+2-2r}}\mathrm{e}^{-\lambda^2_1 T/2}.
\end{equation}
From (\ref{u0first}), (\ref{uve}) and (\ref{entire}) we obtain
\begin{equation}\label{logF}
	\log |F(s)|\leq \frac{T}{2}|\Im(s)-\delta|+\log\pts{(\mathcal{K}+\varepsilon) T^{1 / 2}\frac{\abs{J'_\nu\pts{j_{\nu,1}}}}{\pts{2\kappa_\al}^{1/2}}},\quad s\in\mathbb{C}.
\end{equation}
We recall the following representation theorem, see \cite[p. 56]{koosis}.
\begin{theorem}\label{repre} Let $g(\cdot)$ be an entire function of exponential type and assume that
	$$
	\int_{-\infty}^{\infty} \frac{\log ^{+}|g(x)|}{1+x^{2}} \mathrm{d}x<\infty.
	$$
	Let $\left\{b_{m}\right\}_{m \geq 1}$ be the set of zeros of $g(z)$ in the upper half plane $\Im(z)>0$ (each zero being repeated as many times as its multiplicity). Then,
	$$
	\log |g(z)|=A \Im(z)+\sum_{m=1}^{\infty} \log \left|\frac{z-b_{m}}{z-\bar{b}_{m}}\right|+\frac{\Im(z)}{\pi} \int_{-\infty}^{\infty} \frac{\log |g(s)|}{|s-z|^{2}} \mathrm{d}s,\quad\Im(z)>0,
	$$
	where
	$$
	A=\limsup _{y \rightarrow\infty} \frac{\log |g(i y)|}{y} .
	$$
\end{theorem}
We apply the last result to the function $F$ given in (\ref{entire}). In this case, (\ref{uve}) implies that $A\leq T/2$. Also notice that $\Im\left(a_{k}\right)>0$, $k\geq 1$, to get
\begin{equation}\label{aprep}
	\log \left|F\left(a_{1}\right)\right|\leq\left(\lambda^2_1+\delta\right)\frac{T}{2}+\sum_{k=2}^{\infty} \log \left|\frac{a_{1}-a_{k}}{a_{1}-\bar{a}_{k}}\right|+\frac{\Im\left(a_{1}\right)}{\pi} \int_{-\infty}^{\infty} \frac{\log |F(s)|}{\left|s-a_{1}\right|^{2}} \mathrm{~d}s.
\end{equation}
By using the definition of the constants $a_k$'s we have
\begin{eqnarray}
	\sum_{k=2}^{\infty} \log \left|\frac{a_{1}-a_{k}}{a_{1}-\bar{a}_{k}}\right|&=&\sum_{k=2}^{\infty} \log \left(\frac{\left( j_{\nu, k}\right)^{4}-\left( j_{\nu, 1}\right)^{4}}{2 \delta / \kappa_\alpha^{4}+\left( j_{\nu, 1}\right)^{4}+\left(j_{\nu, k}\right)^{4}}\right)\notag\\
	&\leq& \sum_{k=2}^{\infty} \frac{1}{j_{\nu, k+1}-j_{\nu, k}} \int_{j_{\nu, k}}^{j_{\nu, k+1}} \log \left(\frac{ x^{4}}{2 \delta / \kappa_\alpha^{4}+ x^{4}}\right) \mathrm{d} x \notag\\ 
	&\leq& \frac{1}{\pi} \int_{j_{\nu, 2}}^{\infty} \log \left(\frac{ x^{4}}{2 \delta / \kappa_\alpha^{4}+x^{4}}\right) \mathrm{d} x,\label{apoyo1}\\
	&=& -\frac{j_{\nu,2}}{\pi}\log\pts{\frac{1}{1+2\delta/(\kappa_\al j_{v,2})^4}}-\frac{\sqrt[4]{2\delta}}{\sqrt{2}\pi\kappa_\al}\log\pts{\frac{\kappa^2_\al j^{2}_{\nu,2}-\sqrt[4]{8\delta}\kappa_\al j_{\nu,2}+\sqrt{2\delta}}{\kappa^2_\al j^{2}_{\nu,2}+\sqrt[4]{8\delta}\kappa_\al j_{\nu,2}+\sqrt{2\delta}}}\notag\\
	& & -\frac{\sqrt[4]{8\delta}}{\pi\kappa_\al}\cts{\tan^{-1}\pts{1-\frac{\sqrt{2}\kappa_\al j_{\nu,2}}{\sqrt[4]{2\delta}}}- \tan^{-1}\pts{1+\frac{\sqrt{2}\kappa_\al j_{\nu,2}}{\sqrt[4]{2\delta}}}} -\frac{\sqrt[4]{8\delta}}{\kappa_\al},\notag
\end{eqnarray}
where we have used Lemma \ref{consec} and made the change of variables
$$
\tau=\frac{ \kappa_\alpha}{\sqrt[4]{2 \delta}} x.
$$

From (\ref{logF}) we get the estimate
\begin{equation}\label{apoyo2}
	\frac{\Im\left(a_{1}\right)}{\pi} \int_{-\infty}^{\infty} \frac{\log |F(s)|}{\left|s-a_{1}\right|^{2}} \mathrm{~d} s \leq \frac{T \delta}{2}+\log \left((\mathcal{K}+\varepsilon) T^{1 / 2} \frac{\left|J_{\nu}^{\prime}\left(j_{\nu, 1}\right)\right|}{\pts{2 \kappa_\alpha}^{1/2}}\right).
\end{equation}
From (\ref{ena}), (\ref{aprep}), (\ref{apoyo1}) and (\ref{apoyo2}) we have
\begin{equation}\label{abajo}
	\begin{array}{rl}
		\d\frac{\sqrt[4]{8\delta}}{\kappa_\al}-\frac{j_{\nu,2}}{\pi}\log\pts{1+\frac{2\delta}{\pts{\kappa_\al j_{v,2}}^4}} + \frac{\sqrt[4]{2\delta}}{\sqrt{2}\pi\kappa_\al}\log\pts{\frac{\kappa^2_\al j^{2}_{\nu,2}-\sqrt[4]{8\delta}\kappa_\al j_{\nu,2}+\sqrt{2\delta}}{\kappa^2_\al j^{2}_{\nu,2}+\sqrt[4]{8\delta}\kappa_\al j_{\nu,2}+\sqrt{2\delta}}} & \\[0.5cm]
		&\hspace{-11cm}\d +\hspace{0.2cm}\frac{\sqrt[4]{8\delta}}{\pi\kappa_\al}\tan^{-1}\pts{1-\frac{\sqrt{2}\kappa_\al j_{\nu,2}}{\sqrt[4]{2\delta}}}
		-\frac{\sqrt[4]{8\delta}}{\pi\kappa_\al}\tan^{-1}\pts{1+\frac{\sqrt{2}\kappa_\al j_{\nu,2}}{\sqrt[4]{2\delta}}} 
		-\pts{\lambda_1+\delta}T \\[0.5cm]
		&\hspace{-3cm}\d\leq \log(\mathcal{K}+\varepsilon)+\log h(\alpha,\beta,\mu, T),
	\end{array}
\end{equation}
where
\[
\begin{array}{rcl}
	\d h(\alpha,\beta,\mu, T) & = &\d\frac{\pts{{2T \kappa^{5-4r}_\alpha}}^{1/2}\cts{\pts{\sqrt{\mu(\alpha+\beta)}+\kappa_\al \nu}(1-\ell)+\ell}\left(j_{\nu, 1}\right)^{\nu+2-2r}}{2^{\nu} \Gamma(\nu+1) \left|J_{\nu}^{\prime}\left(j_{\nu, 1}\right)\right|}\\
	& &\\
	&=&\d\frac{\pts{{2T \kappa^{5-4r}_\alpha}}^{1/2}\cts{\pts{\sqrt{\mu(\alpha+\beta)}+\sqrt{\mu(\alpha+\beta)-\mu}}(1-\ell)+\ell}\left(j_{\nu, 1}\right)^{\nu+2-2r}}{2^{\nu} \Gamma(\nu+1) \left|J_{\nu}^{\prime}\left(j_{\nu, 1}\right)\right|}.
\end{array}
\]
The result follows by taking 
\[\delta=2\kappa_\alpha^4 \pts{j_{\nu,2}}^4, \quad \text{and then letting}\quad \varepsilon\rightarrow 0^+.\]
\hfill $\square$


\section{Comments about case $\al+\beta = 1$}\label{sec5}
In this section, we consider the case $\al+\beta= 1$. In the previous sections, we proved the well-posedness of our system (\ref{problem}) for the cases $\alpha + \beta <1$ and $\alpha + \beta >1$. In both cases, we used a Hardy generalized inequality. However, in the case $\al +\beta =  1$, this inequality does not provide any useful information, since $\al + \beta =1$, then $\mu(\al+\beta) = 0$. Hence, since in \cite{GaloLopez2} was studied for the operator (\ref{opeA}) using the Sturm-Liouville theory for this case, we can proceed as in \cite[section 5]{GaloLopez2} or \cite{GaloLopez3} to prove the well-posedness of the systems (\ref{problem}) relative to the operator $\A^2$ given by (\ref{A_square}) and (\ref{A2consts}).\\

From \cite[Section 5]{GaloLopez2} or \cite{Zettl}, if we assume that $0\leq \al<2$, $\beta\in\R$ such that $\beta = 1-\al$, and $\mu < 0$ and consider the differential expression $M$ defined by
\[Mu=-(pu')'+qu\]
where $\d p(x) = x, q(x) = -\mu x^{-1}, w(x) = x^{1-\al}$, $'=\frac{\mathrm{d}}{\mathrm{d}x}$.\\

We can see that,
\begin{equation*}\label{Acond1}
	1/p, q, w \in L_{\text{loc}}(0,1),\quad\text{and}\quad p,w >0\text{ on } (0,1),
\end{equation*} 
thus $Mu$ is defined a.e. for functions $u$ such that $u, pu'\in AC_{\text{loc}}(0,1)$, where $AC_{\text{loc}}(0,1)$ is the space of all locally absolutely continuous functions in $(0,1)$.\\

Notice that the operator $\mathcal{A}$ given in (\ref{opeA}), with $\beta = 1-\al$, can be written as $\A =w^{-1}M$. Now, consider
\begin{equation*}\label{Dmax}
	D_{\max}:=\left\{u\in AC_{\text{loc}}(0,1)\, |\,pu'\in AC_{\text{loc}}(0,1),\, u, \A u\in L^2_{1-\al}(0,1)\right\},\quad\text{and}
\end{equation*}
\[D(\A):=\left\{\begin{aligned}
	\{u\in D_{\max}\,| \lim_{x\rightarrow 0^+}x^{\sqrt{-\mu}}u(x)=0,u(1)=0\} & &\text{if } \sqrt{-\mu}<\kappa_\al, \\
	\{u\in D_{\max}\,| u(1)=0\}& & \text{if } \sqrt{-\mu}\geq\kappa_\al.
\end{aligned}\right.
\]
Recall that the Lagrange form associated with $M$ is given as follows
\[[u,v]:=upv'-vpu',\quad u,v\in D_{\max}.\]

From \cite[Proposition 20]{GaloLopez2} we have the following result
\begin{proposition}
	Let $0\leq \al<2$, $\mu < 0$, and $\nu=\sqrt{-\mu}/\kappa_\al$. Then $\A:D(\A)\subset L^2_{1-\al}(0,1)\rightarrow L^2_{1-\al}(0,1)$ is a self-adjoint operator. Furthermore, when $\beta = 1-\al$ the family given in (\ref{Phik}) is an orthonormal basis for $L^2_{1-\al}(0,1)$ such that
	\begin{equation*}
		\mathcal{A} \Fi_k=\lambda_k \Fi_k, \quad k\geq 1.
	\end{equation*}
\end{proposition}
Then for the operator $\A^{2}$ given by (\ref{A_square}) we have that
\begin{proposition}\label{base1}
	The operator
	\begin{equation}\label{Asquard_1}
		-\A^{2}:D(\A^2)\subset L^{2}_{\beta}(0,1) \rightarrow L^{2}_{\beta}(0,1)
	\end{equation} 
	with dense domain
	
	\begin{equation}\label{domainA21}
		D(\A^2):=\left\{\begin{array}{ll}
			\begin{array}{l}
				\d \left\{u\in D_{\max}\,|\A u\in D_{\max},\, \lim_{x\rightarrow 0^+}x^{\sqrt{-\mu}}u(x)\right.\\ 
				\d\left.\hspace{3cm}=\lim_{x\rightarrow 0^+}x^{\sqrt{-\mu}}\A u(x)=\pts{xu_{x}}_{x}(1) = u(1)=0\right\}
			\end{array}  &\text{if } \sqrt{-\mu}<\kappa_\al, \\[1.2 cm]			
			\{u\in D_{\max}\,| \pts{xu_{x}}_{x}(1) = u(1)=0\}&  \text{if } \sqrt{-\mu}\geq\kappa_\al,
		\end{array}\right.
	\end{equation}
	is a negative self-adjoint operator. Furthermore, the family given by (\ref{Phik})
	is an orthonormal basis for $L^2_\beta(0,1)$ such that it satisfies (\ref{lambda2_k}).
\end{proposition} 
\begin{proof}
	The proof is followed in view of the properties of the operator (\ref{opeA}) and the last proposition.
\end{proof}

Now, when $\al+\beta =  1$ similar to Section \ref{sec2}, we consider convenient interpolation spaces $\HH^s$, $s\in\R$.

\begin{definition}
	Let $T>0$ and $\al,\mu\in \R$ with $0\leq\al<2$, $\mu<0$, and $r\in\{0,1\}$. Let $f \in L^2(0,T)$ and $u_0\in \HH^{-s}$ for some $s > 0$. Then, a weak solution of (\ref{problem}) is a function $u \in C^0([0,T];\HH^{-s})$ such that for every $\tau \in (0,T]$ and for
	every $z^\tau \in \HH^s$ we have
	\begin{equation}\label{weaksolD01}
		\left\langle u(\tau), z^{\tau}\right\rangle_{\mathcal{H}^{-s}, \mathcal{H}^{s}}=\left\langle u_{0}, S(\tau)z^\tau\right\rangle_{\mathcal{H}^{-s}, \mathcal{H}^{s}} - \int_{0}^{\tau} f(t) \mathcal{O}_{1-\sqrt{-\mu}}\pts{\partial_{x}\cts{\A^{1-r} S(\tau-t)z^{\tau}}}\mathrm{d} t,
	\end{equation}
\end{definition}
The last definition is the same in (\ref{weaksol}) when it is considered $\al+\beta = 1$ and $\ell = 0$.

\begin{proposition}\label{continuityD01}
	Let $T>0$, $\al,\mu, r\in \R$ with $0\leq\al<2$ given as in the last definition. Let $f \in L^2(0,T)$ and $u_0\in \HH^{-s}$ such that $s>(\sqrt{-\mu}+\kappa_\al-2r)/2\kappa_\al$. Then, formula (\ref{weaksolD01}) defines for each $\tau \in [0, T ]$ a unique element $u(\tau) \in \HH^{-s}$ that can be written as
	\[
	u(\tau)=S(\tau) u_{0}-\hat{B}_{r}(\tau) f, \quad \tau \in(0, T],\]
	where $B(\tau)$ is the strongly continuous family of bounded operators $\hat{B}_{r}(\tau): L^{2}(0,T) \rightarrow \mathcal{H}^{-s}$ given by
	\[\left\langle \hat{B}_{r}(\tau) f, z^{\tau}\right\rangle_{\mathcal{H}^{-s}, \mathcal{H}^{s}}=\int_{0}^{\tau} f(t) \mathcal{O}_{1-\sqrt{-\mu}}\pts{\partial_{x}\cts{\A^{1-r} S(\tau-t)z^{\tau}}}\mathrm{d} t, \quad \text{for all  }z^{\tau} \in \mathcal{H}^{s} .\]
	Furthermore, the unique weak solution $u$ on $[0, T]$ to (\ref{problem}) (in the sense of (\ref{weaksolD01})) belongs to ${C}^{0}\left([0, T] ; \mathcal{H}^{-s}\right)$ and fulfills
	\[
	\|u\|_{L^{\infty}\left([0, T] ; \mathcal{H}^{-s}\right)} \leq C\left(\left\|u_{0}\right\|_{\mathcal{H}^{-s}}+\|f\|_{L^{2}(0, T)}\right).
	\]
\end{proposition}

We are in the same position as in the case $\al+\beta < 1$, so we can follow the same steps to get following result,
\begin{theorem}\label{Teo2}
	Let $T>0$, and $\al,\, \beta,\, \mu\in \R$ with $0\leq\al<2$, $\al+\beta = 1$, $\mu < 0$ and $r\in\{0,1\}$. The next statements hold.
	\begin{enumerate}
		\item \textbf{Existence of a control.} For any $u_0\in L^2_{1-\alpha}(0,1)$ there exists a control $f \in L^2(0, T )$ such that the solution $u$ to (\ref{problem}) satisfies $u(\cdot,T ) \equiv 0$.
		\item \textbf{Upper bound of the cost.} There exists a constant $c>0$ such that for every $\delta\in (0,1)$ we have that when $\gamma$ is given by (\ref{gamma}), then 
		\[\mathcal{K}(T,\alpha,\mu,r)\leq \frac{c M(T,\alpha,\nu,\delta)T^{1/2}}{\kappa_\al^{9/2-2r}\sqrt{-\mu}}
		\exp\pts{-\frac{T}{2}\kappa_\alpha^4 j_{\nu,1}^4}.\]
		where $\kappa_\al$ is given by (\ref{Nu}) and $j_{\nu,1}$ is the first positive zero of the Bessel function $J_\nu$ (defined in the Appendix), and
		\[
		\begin{array}{rcl}
			M(T,\alpha,\nu,\delta)& = &\d\pts{1+\frac{1}{(1-\delta)\kappa_\alpha^2 T}}\cts{\exp\pts{\frac{\sqrt{2+\sqrt{2}}}{\sqrt{2}\kappa_\alpha}}+\frac{1}{\delta^3}\exp\pts{\frac{3\sqrt{2+\sqrt{2}}}{(1-\delta)\kappa_\alpha^2 T}}}\\
			&&\\
			&&\d\times\exp\pts{-\frac{(1-\delta)^{3/2}T^{3/2}}{8\sqrt{2+\sqrt{2}}(1+T)^{1/2}}\kappa_\alpha^5 j_{\nu,1}^4}.
		\end{array}
		\]
		
		\item \textbf{Lower bound of the cost.} There exists a constant $c>0$ such that
		\[\frac{c2^{\nu} \Gamma(\nu+1) \left|J_{\nu}^{\prime}\left(j_{\nu, 1}\right)\right|\exp{\left(2\pts{1-\frac{\log 5}{\pi}-\frac{\tan^{-1} 2}{\pi}}j_{\nu,2}\right)}}{\pts{{2T \kappa^{5-4r}_\alpha}}^{1/2}\left(j_{\nu, 1}\right)^{\nu+2-2r}\sqrt{-\mu}}\exp\pts{-\pts{j_{\nu,1}^4+2j_{\nu,2}^4}\kappa_\alpha^4 T}\leq \mathcal{K}(T,\alpha,\mu, r),\]			
	\end{enumerate}
	where $j_{\nu,2}$ is the second positive zero of the Bessel function $J_\nu$.
	
\end{theorem}


\appendix
\section{Bessel functions}\label{appen}

We introduce the Bessel function of the first kind $J_{\nu}$ as follows
\begin{equation}\label{bessel}
	J_{\nu}(x)=\sum_{m \geq 0} \frac{(-1)^{m}}{m ! \Gamma(m+\nu+1)}\left(\frac{x}{2}\right)^{2 m+\nu}, \quad x \geq 0,
\end{equation}
where $\Gamma(\cdot)$ is the Gamma function. In particular, for $\nu>-1$ and $0<x \leq \sqrt{\nu+1}$, from (\ref{bessel}) we have (see \cite[9.1.7, p. 360]{abram})
\begin{equation}\label{asincero}
	J_{\nu}(x) \sim \frac{1}{\Gamma(\nu+1)}\left(\frac{x}{2}\right)^{\nu} \quad \text { as } \quad x \rightarrow 0^{+} .
\end{equation}
A Bessel function $J_\nu$ of the first kind solves the differential equation
\begin{equation}\label{Besselode}
	x^2y''+xy'+(x^2-\nu^2)y=0.
\end{equation}
Bessel functions of the first kind satisfy the recurrence formula (\cite{abram}, 9.1.27):
\begin{equation}\label{recur}
	x J_{\nu}^{\prime}(x)-\nu J_{\nu}(x)=-x J_{\nu+1}(x).
\end{equation}

Recall the asymptotic behavior of the Bessel function $J_{\nu}$ for large $x$, see \cite[Lem. 7.2, p. 129]{komo}.
\begin{lemma}\label{asimxinf}
	For any $\nu \in \mathbb{R}$
	$$
	J_{\nu}(x)=\sqrt{\frac{2}{\pi x}}\left\{\cos \left(x-\frac{\nu \pi}{2}-\frac{\pi}{4}\right)+\mathcal{O}\left(\frac{1}{x}\right)\right\} \quad \text { as } \quad x \rightarrow \infty.
	$$
\end{lemma}

For $\nu >0$ the Bessel function $J_{\nu}$ has an infinite number of real zeros $(j_{\nu,k})$ all of which are simple, with the possible exception of $x=0$, see \cite[9.5.2, p. 370]{abram}, satisfying
\begin{equation}\label{cerocrec} 
	\nu<j_{\nu, 1}<j_{\nu, 2}<\ldots.
\end{equation}
Moreover, from \cite[Chap. XV, p. 486, eq. (5)]{Watson}, the first zero $j_{\nu,1}$ satisfies
$$
\sqrt{\nu(\nu+2)} < j_{\nu,1} < \sqrt{2(\nu+1)(\nu+3)},
$$
then, if $\nu > 1$, we have an upper bound of $j_{v,1}$ given by 
\begin{equation}\label{BoundFirstZero}
	\nu < j_{\nu,1} < \sqrt{15}\,\nu+\frac{1}{\sqrt{\nu}}.
\end{equation}

We recall Stirling's formula, see \cite[6.1.39, p. 257]{abram},
\begin{equation}\label{stirling}
	\Gamma(\nu+1)\sim \sqrt{2\pi\nu}\left(\frac{\nu}{\mathrm{e}}\right)^{\nu}\quad\text{as } \nu\to\infty.
\end{equation}
In \cite[9.3.1, p. 365]{abram}, we have
\begin{equation}\label{Asymtotic}
	J_{\nu}(x)\sim \frac{1}{\sqrt{2\pi\nu}}\left(\frac{\mathrm{e}x}{2\nu}\right)^{\nu}\quad\text{as } \nu\to\infty.
\end{equation}

For $\nu >-1$ the Bessel function $J_{\nu}$ has an infinite number of real zeros $0<j_{\nu, 1}<j_{\nu, 2}<\ldots$, all of which are simple, with the possible exception of $x=0$. In \cite[Proposition 7.8]{komo} we can find the next information about the location of the zeros of the Bessel functions $J_{\nu}$:
\begin{lemma}\label{consec}Let $\nu \geq 0$.\\
	1. The difference sequence $\left(j_{\nu, k+1}-j_{\nu, k}\right)_{k}$ converges to $\pi$ as $k \rightarrow\infty$.\\
	2. The sequence $\left(j_{\nu, k+1}-j_{\nu, k}\right)_{k}$ is strictly decreasing if $|\nu|>\frac{1}{2}$, strictly increasing if $|\nu|<\frac{1}{2}$, and constant if $|\nu|=\frac{1}{2}$.\\
\end{lemma}

For $\nu \geq 0$ fixed, we consider the next asymptotic expansion of the zeros of the Bessel function $J_{\nu}$, see \cite[Section 15.53]{Watson},
\begin{equation}\label{asint}
	j_{\nu, k}=\left(k+\frac{\nu}{2}-\frac{1}{4}\right) \pi-\frac{4 \nu^{2}-1}{8\left(k+\frac{\nu}{2}-\frac{1}{4}\right) \pi}+O\left(\frac{1}{k^{3}}\right), \quad \text { as } k \rightarrow\infty.
\end{equation}

In particular we have
\begin{equation}\label{below}
	\begin{aligned}
		&j_{\nu, k} \geq\left(k-\frac{1}{4}\right) \pi \quad \text { for } \nu \in\left[0, 1/2\right], \\
		&j_{\nu, k} \geq\left(k-\frac{1}{8}\right) \pi \quad \text { for } \nu \in\left[1/2,\infty\right).
	\end{aligned}
\end{equation}

\begin{lemma}\label{reduce} For any $\nu \geq 0$ and any $k\geq 1$ we have
	$$
	\sqrt{j_{\nu, k}}\left|J_{\nu}^{\prime}\left(j_{\nu, k}\right)\right|=\sqrt{\frac{2}{\pi}}+O\left(\frac{1}{j_{\nu, k}}\right)\quad \text{as}\quad k \rightarrow \infty.
	$$
\end{lemma}
The proof of this result follows by using  (\ref{asincero}) and the recurrence formula (\ref{recur}).

\begin{lemma} Let $\ell=\ell(\al,\beta)$, $\gamma=\gamma(\alpha,\beta,\mu)$ and $\nu=\nu(\alpha,\beta,\mu)$ given in (\ref{index_l}), (\ref{gamma}) and (\ref{Nu}) respectively, then the $\al+\beta+\gamma-\ell$-generalized limit of $d^{1-\ell}\Fi_k/dx^{1-\ell}=\Fi^{(1-\ell)}_k$ at $x=0$ is finite for all $k\geq 1$, and
	\begin{equation}\label{limAux}
		\OO_{\al+\beta+\gamma-\ell}(\Fi^{(1-\ell)}_k)= \frac{(2\kappa_\al)^{1/2}\cts{\pts{\sqrt{\mu(\alpha+\beta)}+\kappa_\al \nu}(1-\ell)+\ell}\pts{j_{\nu,k}}^{\nu}}{2^{\nu}\Gamma(\nu+1)\abs{J'_{\nu}\pts{j_{\nu, k}}}}, \quad k\geq1.
	\end{equation}
\end{lemma}
\begin{proof}
	If $\ell=0$, from the recurrence formula (\ref{recur}) we obtain that
	$$
	\Fi_{k}'(x) =\frac{(2\kappa_\al)^{1/2}}{\abs{J'_\nu\pts{j_{\nu,k}}}}\cts{ \pts{\frac{1-\al-\beta}{2}+\kappa_\al\nu}x^{-(1+\al+\beta)/2}J_{\nu}\pts{j_{\nu,k}x^{\kappa_\al}}-\kappa_\al j_{\nu,k}x^{-(1+\al+\beta)/2+\kappa_\al}J_{\nu+1}\pts{j_{\nu,k}x^{\kappa_\al}}},
	$$
	and using  (\ref{asincero}), we obtain that
	$$
	\OO_{\al+\beta+\gamma}(\Fi'_k)= \frac{(2\kappa_\al)^{1/2}\pts{\sqrt{\mu(\alpha+\beta)}+\kappa_\al \nu}\pts{j_{\nu,k}}^{\nu}}{2^{\nu}\Gamma(\nu+1)\abs{J'_{\nu}\pts{j_{\nu, k}}}}.
	$$
	For case $\ell=1$, we only use (\ref{asincero}) to obtain
	$$
	\mathcal{O}_{\alpha+\beta+\gamma-1}(\Phi_k)= \frac{(2\kappa_\alpha)^{1/2}(j_{\nu,k})^{\nu}}{2^{\nu}\Gamma(\nu+1)|J'_{\nu}(j_{\nu, k})|}.
	$$ 
\end{proof}

Since (\ref{bessel}) considers the modified Bessel function
\begin{equation}\label{besselModified}
	I_{\nu}(x) := i^{-\nu}J_{\nu}(ix) = \sum_{m \geq 0} \frac{1}{m ! \Gamma(m+\nu+1)}\left(\frac{x}{2}\right)^{2 m+\nu}, \quad x \geq 0, 
\end{equation}
where $\nu$ is given by (\ref{Nu}).

From (\ref{besselModified}) we have the following lower bound to the modified Bessel function,
\begin{equation}\label{lowerBM}
	I_{\nu}(x)\geq \dfrac{x^{\nu}}{2^{\nu}\Gamma(\nu+1)}.
\end{equation}






\end{document}